\documentclass{article}

\usepackage{graphicx}
\usepackage{mathrsfs}
\usepackage{CJK}
\usepackage{amsmath}
\usepackage{amsthm}
\usepackage{amsfonts}
\usepackage{indentfirst} 
\usepackage{latexsym,bm}
\usepackage{mathrsfs}
\usepackage{cite}
\usepackage{fancyhdr}
\usepackage{fancybox}
\usepackage[noend]{algpseudocode}
\usepackage[ruled,linesnumbered]{algorithm2e}
\usepackage{geometry}
\newtheorem{definition}{Definition}

\newtheorem{theorem}{Theorem}
\newtheorem{lemma}{Lemma}
\usepackage{titling}
\usepackage{color}
\usepackage{hyperref}
\usepackage{footmisc}

\title{Generalized Lagrangian Neural Networks}

\begin{document} 

\begin{CJK}{UTF8}{gbsn} 
\date{}

\maketitle
\hspace*{\fill}\\ 
\begin{center}
\textbf{Shanshan Xiao} \footnote{LSEC, ICMSEC, Academy of Mathematics and systems Science, Chinese Academy of Sciences, Beijing 100190, China \label{lsec}} \footnote{School of Mathematical Sciences, University of Chinese Academy of Sciences, Beijing 100049, China \label{ucas}} $\cdot$ \textbf{Jiawei Zhang}\footref{lsec} \footref{ucas} $\cdot$ \textbf{Yifa Tang}$^{*}$ \footref{lsec} \footref{ucas}
\end{center}

\hspace*{\fill}\\\\\\

\begin{abstract}
    Incorporating neural networks for the solution of Ordinary Differential Equations (ODEs) represents a pivotal research direction within computational mathematics. Within neural network architectures, the integration of the intrinsic structure of ODEs offers advantages such as enhanced predictive capabilities and reduced data utilization. Among these structural ODE forms, the Lagrangian representation stands out due to its significant physical underpinnings. Building upon this framework, Bhattoo introduced the concept of Lagrangian Neural Networks (LNNs). Then in this article, we introduce a groundbreaking extension (Genralized Lagrangian Neural Networks) to Lagrangian Neural Networks (LNNs), innovatively tailoring them for non-conserative systems. By leveraging the foundational importance of the Lagrangian within Lagrange's equations, we formulate the model based on the generalized Lagrange's equation. This modification not only enhances prediction accuracy but also guarantees Lagrangian representation in non-conservative systems. Furthermore, we perform various experiments, encompassing 1-dimensional and 2-dimensional examples, along with an examination of the impact of network parameters, which proved the superiority of Genralized Lagrangian Neural Networks(GLNNs).
\end{abstract}


\setlength{\parindent}{2em}

\section{Introduction}

Machine learning has found significant applications in the field of mathematics, revolutionizing traditional approaches to problem-solving and analysis. With its ability to automatically learn patterns and make predictions from data, machine learning techniques have been successfully applied to various mathematical tasks. One significant direction is the utilization of machine learning to address mathematical problems associated with differential equations and dynamical systems\cite{cranmer2020lagrangian, bersani2021lagrangian,zhong2020unsupervised,chen2018neural,chen2019symplectic, jin2020sympnets,tong2021symplectic,greydanus2019hamiltonian,yu2021onsagernet}.

In the existing methods for solving dynamical system problems using machine learning, they can be mainly categorized into two types. One type is unstructured methods\cite{chen2018neural,greydanus2019hamiltonian, gomez2017reversible, mackay2018reversible, quade2016prediction, chan1999galerkin, sherstinsky2020fundamentals}, which do not consider the physical or mathematical structure of the equations or dynamical systems. When employing this approach, the focus is typically on minimizing the model error or complexity. For example, symbolic regression\cite{quade2016prediction} is a regression method where the control equation under study is treated as an unknown target, and this unknown target is regarded as a function of the state variable data and its time derivative, using certain sparse approximations; Galerkin-closure methods\cite{chan1999galerkin}] aim to approximate the unresolved scales of turbulence by using a closure model based on a truncated set of resolved scales; furthermore, unstructured neural networks such as LSTM\cite{sherstinsky2020fundamentals} can be used to directly predict the solutions of differential equations or the phase flow of dynamical systems.

In contrast to unstructured methods, considering the inherent properties of the system and utilizing neural networks with specific structures can significantly enhance the predictive performance on specific systems. By incorporating the knowledge of system properties, such as conservation laws, symmetries, or known mathematical structures, into the design of neural networks, it becomes possible to improve the accuracy and effectiveness of predictions for the targeted system\cite{cranmer2020lagrangian, zhong2020unsupervised, chen2019symplectic, jin2020sympnets, tong2021symplectic, greydanus2019hamiltonian, lutter2019deep, sanchez2019hamiltonian, toth2019hamiltonian}. Hamiltonian Neural Network (HNN)\cite{greydanus2019hamiltonian}, is a typical example of a structured neural network that takes into account the geometric structure of Hamiltonian systems. It utilizes neural networks to approximate the Hamiltonian of the system, thereby achieving improved predictive performance. OnsagerNets\cite{yu2021onsagernet}, as a systematic method that can overcome the aforementioned limitations, are base on a highly general extension of the Onsager principle for dissipative dynamics.
Lagrangian Neural Networks (LNN)\cite{cranmer2020lagrangian} are a class of neural networks specifically designed to parameterize arbitrary Lagrangians through their network architecture. Unlike traditional approaches, LNNs do not impose restrictions on the functional form of the learned energies, allowing them to produce models that conserve energy. 

LNNs have demonstrated exceptional performance in many tasks. However, LNNs are limited to Lagrangian systems that adhere to the principle of energy conservation. However, in practical problems, we often encounter non-conservative systems. In the study of Euler-Lagrange equations\cite{santilli1982foundations, santilli2013foundations}, it is known that non-conservative systems can be formulated in the form of generalized Euler-Lagrange equations\cite{bersani2021lagrangian}, which can be expressed as follows: 
$$\frac{d}{dt}\left( \frac{\partial \mathscr{L}}{\partial \Dot{q_{k}}}\right) - \frac{\partial \mathscr{L}}{\partial q_{k}} = F_{k}.$$

In this article, we extend the scope of LNNs by constructing neural networks specifically designed for non-conservative systems. Our main inspiration stems from the special physical significance of the Lagrangian in Lagrange's equations. Since the Lagrangian can be understood as the difference between kinetic energy and potential energy, authors replace the baseline model with a fitted Lagrangian in Lagrangian Neural Networks (LNNs). The physical interpretation of the Lagrangian and the structure of Lagrange's equations contribute to higher prediction accuracy and enable the energy of non-conservative systems to fluctuate within a narrow range (as opposed to energy dissipation in the baseline model). Building upon this, we consider the generalized Lagrange's equations, which incorporate non-conservative terms, akin to the conventional Lagrange's equations. In the generalized Lagrange's equations, the Lagrangian $\mathscr{L}$ holds the same physical interpretation. Thus, our objective is to construct a neural network suitable for non-conservative systems using the framework of generalized Lagrange's equations. To achieve this, in Chapter 2, we first present the mathematical and physical properties of Lagrange's equations, along with methods to express most physically motivated non-conservative systems in the form of generalized Lagrange's equations. Consequently, we can regard the majority of physical non-conservative systems as governed by a generalized Lagrange's equation.

Building upon this theoretical foundation, we introduce the Generalized Lagrangian Neural Networks (GLNNs). We conceptualize non-conservative systems as dynamic systems governed by generalized Lagrange's equations, and employ neural networks that preserve the structural integrity of these generalized equations to model such dynamics. Given the physical interpretations embedded within the generalized Lagrange's equations, GLNNs demonstrate enhanced predictive capabilities in specific scenarios.

The structure of our paper unfolds as follows: In Chapter 2, we elucidate the fundamental concepts of Lagrangian equations, emphasizing that the majority of physical systems can be represented in the form of generalized Lagrangian equations. Additionally, we present a methodology for formulating these generalized Lagrangian equations. Subsequently, Chapter 3 delves into the construction and comparison of GLNNs against a baseline model. Our numerical experiments, focusing on predictive analyses of selected one-dimensional and two-dimensional physical models while varying network hyperparameters, are detailed in Chapter 4. Finally, Chapter 5 encapsulates our findings, highlighting both the merits of GLNNs and certain limitations observed during training.

\section{Preliminaries}

In this chapter, we will introduce Lagrangian systems, generalized Lagrangian systems, and their associated mathematical and physical properties.

\begin{definition}
We call 
\begin{equation}
    \frac{d}{dt}\left( \frac{\partial \mathscr{L}}{\partial \Dot{q_{k}}}\right) - \frac{\partial \mathscr{L}}{\partial q_{k}} = 0
    \label{el1}
\end{equation}
the Lagrangian description of mechanical system, where $\mathscr{L}$ is the Lagrangian, $q = (q_{1}, \cdots, q_{N})$ the generalized coordinates of the point. 
\label{lagrangian}
\end{definition}

If we denote the kinetic energy of system by $T$, and $U$ the potential energy, then we can suppose that the total energy is $T+U$, and the Lagrangian is $\mathscr{L} = T-U$. The Lagrangian systems presented in Definition \ref{lagrangian} are assumed to be unaffected by external forces, thereby exhibiting inherent energy conservation properties. However, in practical scenarios, systems often experience the influence of external forces, such as frictional forces or artificially applied forces. In response to this, a generalized form of the Lagrange equation, similar in structure to that described in Definition \ref{lagrangian}, can be derived as follows\cite{bersani2021lagrangian}:
\begin{equation}
    \frac{d}{dt}\left( \frac{\partial \mathscr{L}}{\partial \Dot{q_{k}}}\right) - \frac{\partial \mathscr{L}}{\partial q_{k}} = F_{k}
    \label{el2}
\end{equation}
here $F_{k}$ represent non-conservative forces.

In certain special cases, we can express the dissipative external force $F_{k}$ in a specific form, such as $ F_{k} = -a_{k}\Dot{q_{k}}^{n} $, for example:
\begin{itemize}
    \item If the external force corresponds to frictional force, the terms $a_{k}$ represent the coefficients of friction, and it is assumed that $n = 0$.
    \item If the external force corresponds to viscous force, then the terms $a_{k}$ represent the viscosity coefficints, and $n \geq 1.$ One of the most common examples is that of a damped harmonic motion, where the external force is given by $F = -a\Dot{q}$, and the system can be described as $\Ddot{q} + a\Dot{q} + k^{2}q = 0$.
\end{itemize}

After providing the definition of a Lagrangian system, we will reference several theorems to illustrate the types of systems that can be represented in the form of (generalized) Lagrange's equations. Furthermore, we will present some expressions of Lagrange's equations for specific systems.

\begin{theorem}[Fundamental analytic theorem for configuration space formulations]\cite{santilli1982foundations}

A necessary and sufficient condition for a local, holonomic, generally nonconservative Newtonian system in the fundamental form
\begin{equation}
    A_{ki}(t,q,\Dot{q})\Ddot{q}^{i} + B_{k}(t,q,\Dot{q}) = 0, k = 1,2,\cdots,n,
\end{equation}
which is well defined, of class $\mathscr{C}^{2}$, and regular in a star-shaped region $\mathbb{R}^{*2n+1}$ of the variables$(t,q,\Dot{q})$, to admit an ordered direct analytic representation in terms of the conventional Lagrange's equation in $\mathbb{R}^{*2n+1}$,
\begin{equation}
    \frac{d}{dt}\frac{\partial L}{\partial \Dot{q}^{k}} - \frac{\partial L}{\partial q^{k}} \equiv   A_{ki}\Ddot{q}^{i} + B_{k},
\end{equation}
is that the system of equations of motion is self-adjoint in $\mathbb{R}^{*2n+1}$.
\label{nl}
\end{theorem}

\begin{theorem}[A method to construct Lagrangian]\cite{santilli1982foundations}

A Lagrangian for the ordered direct analytic representation of local, holonomic, generally nonconservative Newtonian systems that are well defined of class $\mathscr{C}^{2}$,  regular  and self-adjoint in a star-shaped region $\mathbb{R}^{*2n+1}$ of points $(t,q,\Dot{q})$,
\begin{equation}
    A_{ki}(t,q,\Dot{q})\Ddot{q}^{i} + B_{k}(t,q,\Dot{q}) = 0, k = 1,2,\cdots,n,
\end{equation}
is given by
\begin{equation}
    L = K(t,q,\Dot{q}) + D_{k}(t,q)\Dot{q}^{k} + C(t,q)
\end{equation}
where functions $K,D_{k}$ and $C$ are a solution of partial differential equations
\begin{equation}
\begin{split}
     & \frac{\partial^{2} K}{\partial\Dot{q}^{k_{1}}\partial \Dot{q}^{k_{2}}} = A_{k_{1}k_{2}}(t,q,\Dot{q}),\\
     \frac{\partial D_{k_{1}}}{\partial q^{k_{2}}} - \frac{\partial D_{k_{2}}}{\partial q^{k_{1}}} & = \frac{1}{2} \left( \frac{\partial B_{k_{1}}}{\partial q^{k_{2}}} - \frac{\partial B_{k_{2}}}{\partial q^{k_{1}}} \right) + \left( \frac{\partial^{2} K}{\partial q^{k_{1}}\partial \Dot{q}^{k_{2}}} - \frac{\partial^{2} K}{\partial\Dot{q}^{k_{1}}\partial q^{k_{2}}}\right) \\
     & \equiv Z_{k_{1}k_{2}}(t,q), \\
     \frac{\partial C}{\partial q^{k_{1}}} & = \frac{\partial D_{k_{1}}}{\partial t} - B_{k_{1}} - \frac{\partial K}{\partial q^{k_{1}}} + \frac{\partial^{2} K}{\partial\Dot{q}^{k_{1}}\partial t} \\
     & + \left[ \frac{\partial^{2} K}{\partial q^{k_{1}}\partial \Dot{q}^{k_{2}}} +  \frac{1}{2} \left( \frac{\partial B_{k_{1}}}{\partial q^{k_{2}}} - \frac{\partial B_{k_{2}}}{\partial q^{k_{1}}} \right) \right]\Dot{q}^{k_{2}}\\
     & \equiv W_{k}(t,q),
\end{split}
\end{equation}
given by
\begin{equation}
\begin{split}
K(t,q,\Dot{q}) = \Dot{q^{k_{1}}} & \int_{0}^{1}d\tau' \left\{ \left[ \int_{0}^{1}d\tau A_{k_{1}k_{2}}(t,q,\tau\Dot{q}) \right] \right\}(t,q,\tau'\Dot{q}),\\
& D_{k_{1}} = \left[ \int_{0}^{1}d\tau \tau Z_{k_{1}k_{2}}(t,\tau q) \right]q^{k_{2}},\\
& C = \left[ \int_{0}^{1}d\tau W_{k}(t, \tau q ) \right]q^{k}.
\end{split}
\end{equation}
\label{nl2}
\end{theorem}

Theorems \ref{nl} and \ref{nl2} provide conditions under which Newtonian systems, a class of dynamic systems commonly encountered in realistic models, can be expressed in the form of Lagrange's equations. These theorems also provide the transformation formulas that hold when these conditions are satisfied.

Thus, we have established that self-adjoint Newtonian systems can be expressed in the form of Lagrange's equations, as shown in Equation (\ref{el1}), with a general procedure for obtaining the expressions. The next question is, if we consider a non-conservative system, can it be represented by the generalized form of Lagrange's equations as shown in Equation (\ref{el2})? If so, how can we derive the expression of the generalized Lagrange's equations for a non-conservative system?

In order to achieve better energy prediction performance for the system, we employ the physically meaningful quantities of $T$ (kinetic energy) and $U$ (potential energy) to determine the Lagrangian term $\mathscr{L} = T- U$ in the generalized Lagrangian representation of non-conservative systems. Additionally, the non-conservative force can be obtained using equations:
$$ F_{k} = \frac{d}{dt}\left( \frac{\partial \mathscr{L}}{\partial \Dot{q_{k}}}\right) - \frac{\partial \mathscr{L}}{\partial q_{k}}.$$

Thus, by employing the aforementioned approach, we can obtain the generalized Lagrangian representation for non-conservative systems. In other words, we can view the equations of non-conservative systems as generalized Lagrange's equations.

\begin{lemma}
For most physically motivated non-conservative system, its generalized Lagrangian representation can be obtained using the following approach:
\begin{equation}
\begin{split}
    & \mathscr{L} = T - U \\
    &  F_{k} = \frac{d}{dt}\left( \frac{\partial \mathscr{L}}{\partial \Dot{q_{k}}}\right) - \frac{\partial \mathscr{L}}{\partial q_{k}}.
\end{split}
\end{equation}
here $T$ and $U$ can be calculated by physical formula.
\label{generalL}
\end{lemma}

It is worth noting that, for a given system, the expression form of the generalized Lagrange's equation is not unique. Moreover, it has been proven that different generalized Lagrange's equations representing the same system exhibit certain mathematical connections.

Building upon the theoretical foundations described above, we are able to consider most physically meaningful non-conservative system models as determined by generalized Lagrange's equations. Consequently, we can establish non-conservative Lagrangian neural networks based on this framework.

\section{Theory of generalized Lagrangian neural networks}

Considering a non-conservative system $S$ established from a physical model, given an observed dataset $T = \{(x_{i}, x_{i+1})\}$ obtained from system $S$, our objective is to construct a neural network capable of learning system $S$ by training on the available dataset $T$. The goal is to enable the neural network to predict the phase flow of system $S$ from arbitrary initial points.

As per Lemma \ref{generalL}, we understand that system $S$ can be viewed as governed by the following generalized Lagrange's equations:
$$\frac{d}{dt}\left( \frac{\partial \mathscr{L}}{\partial \Dot{q_{k}}}\right) - \frac{\partial \mathscr{L}}{\partial q_{k}} = F_{k}.$$
In the context of generalized LNNs, we treat both $\mathscr{L}$ and the non-conservative term $F$ as black boxes, i.e., we employ neural networks to learn these two functions. In order to proceed with obtaining predictive results, we apply the chain rule to expand the generalized Lagrange's equations as follows:
\begin{equation}
    (\nabla_{\Dot{q}} \nabla_{\Dot{q}}^{\top}\mathscr{L}) \Ddot{q} + (\nabla_{q} \nabla_{\Dot{q}}^{\top}\mathscr{L})\Dot{q} = \nabla_{q}\mathscr{L} + F.
\end{equation}
Here, the notation $\nabla$ is consistent with the convention used in LNNs, i.e., $(\nabla_{\Dot{q}})_{i} \equiv \frac{\partial}{\partial \Dot{q_{i}}}$. By performing matrix transformations on the above equation, we obtain:
\begin{equation}
    \Ddot{q} = (\nabla_{\Dot{q}} \nabla_{\Dot{q}}^{\top}\mathscr{L})^{-1}[\nabla_{q}\mathscr{L} - (\nabla_{q} \nabla_{\Dot{q}}^{\top}\mathscr{L})\Dot{q} + F ]
    \label{ddotq}
\end{equation}

\begin{figure}[htbp]
\centering
\includegraphics[width=1.0\textwidth]{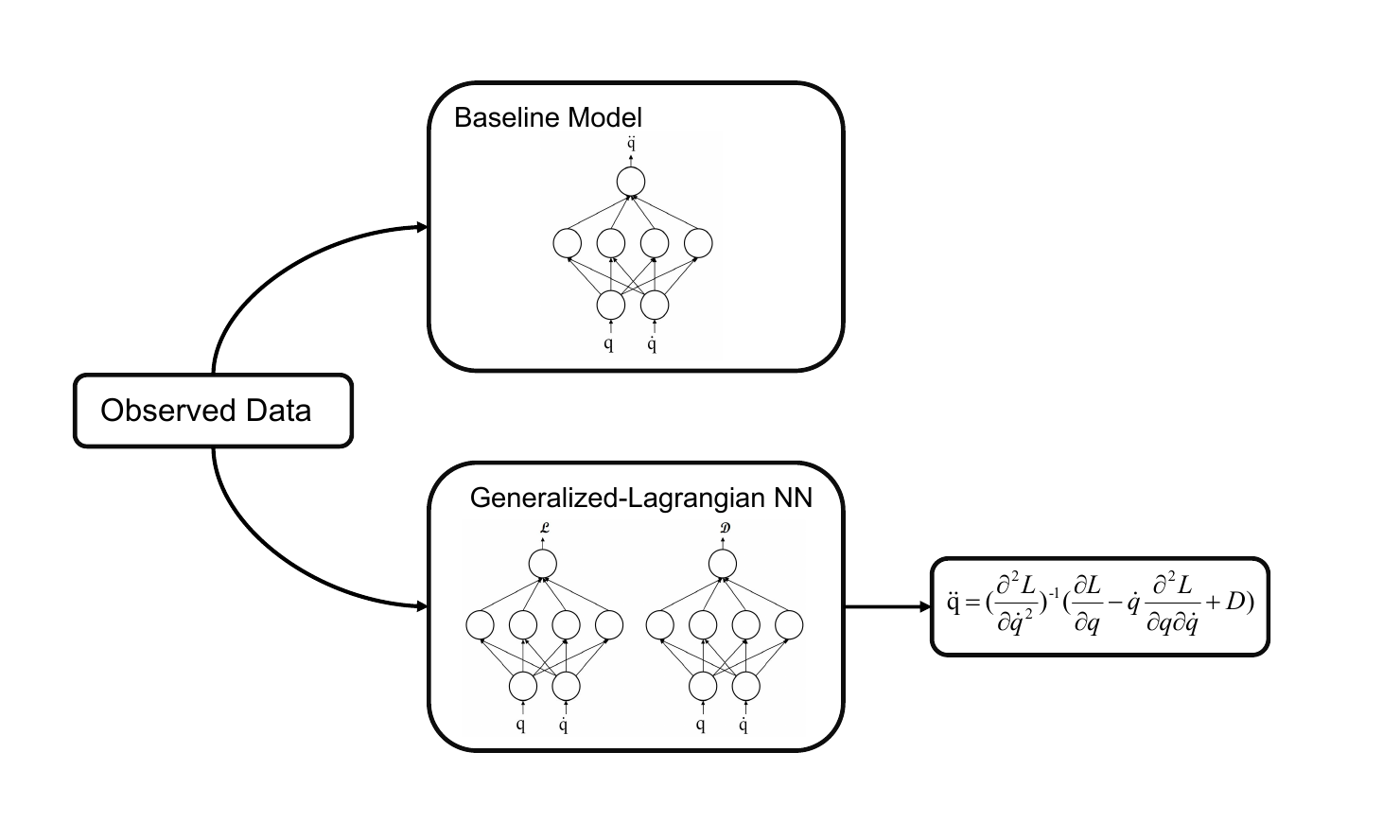}
\caption{Architecture of models}
\label{comparernn}
\end{figure}

Consequently, for the coordinate $x_{t} = (q_{t},\Dot{q_{t}})$ at time $t$, we obtain a method for computing $\Ddot{q_{t}}$ from generalizedd LNNs, allowing us to calculate the phase flow of the system using any numerical scheme.
\\ \hspace*{\fill} \\
\textbf{Loss function}

Here, we provide two possible loss functions, allowing for the selection of the most appropriate one for training based on different input data. When it is possible to obtain or compute the $\ddot{q}$ value at a certain time from raw data, the loss function is chosen as:
$$\mathcal{L} = \frac{1}{|S|} \sum_{x_{t}\in S} \Vert \Tilde{\Ddot{q}} - \Ddot{q} \Vert ^{2},$$
where $x_{t}$ is the coordinate at time $t$, and $\widetilde{\ddot{q}}$ is computed using Equation \ref{ddotq}.
Otherwise, the loss function
$$\mathcal{L} = \frac{1}{|S|} \sum_{x_{t}\in S} \Vert \Tilde{x_{t}} - x_{t} \Vert^{2}$$
is used, where $\widetilde{x_{t}}$ represents the prediction of coordinates by our network, $S$ is the sample data.

\section{Learning a non-conservative system from data}

\subsection{Selection of baseline model}

The examples selected in our paper can be regarded as second-order equations, specifically 
$\Ddot{q} = F(\Dot{q},q,t)$. Therefore, we consider choosing from the following two baseline models:
\begin{itemize}
    \item The first choice is Neural ODE, where we treat the system as determined by the 
    equation $\Dot{q} = f(q,t)$, and we use a neural network to approximate the function $f$. The loss function is $$\mathscr{L} = \frac{1}{|S|} \sum_{x_{k}\in S}\Vert \Tilde{x_{k}} - x_{k} \Vert^{2},$$
    here $\widetilde{x_{t}}$ represents the prediction of coordinates.
    \item Another choice for the baseline model is based on the form of our example equations, which can all be seen as $\Ddot{q} = F(\Dot{q},q,t)$. In this case, we consider Baseline 2: using a neural network to approximate the function $F$, with the same loss function as Neural ODE.
\end{itemize}

Although both of these baseline models can achieve decent prediction performance, due to the shared characteristics of the selected systems, the second baseline model may exhibit superior performance. Consequently, we choose it as the baseline model for comparison.

\begin{table}[htbp]
\caption{Training parameters}
\resizebox{\textwidth}{!}{
\begin{tabular}{|l|llll|}
\hline
Task                                   & Model                                                          & Learning rate                                         & Batch size                                          & Epochs                                              \\ \hline
Damped harmonic motion                 & \begin{tabular}[c]{@{}l@{}}Baseline model\\ GLNNs\end{tabular} & \begin{tabular}[c]{@{}l@{}}0.001\\ 0.001\end{tabular} & \begin{tabular}[c]{@{}l@{}}1000\\ 1000\end{tabular} & \begin{tabular}[c]{@{}l@{}}300\\ 300\end{tabular} \\ \hline
Compound double pendulum with friction & \begin{tabular}[c]{@{}l@{}}Baseline model\\  GLNNs\end{tabular} & \begin{tabular}[c]{@{}l@{}}0.001\\ 0.001\end{tabular} & \begin{tabular}[c]{@{}l@{}}1000\\ 1000\end{tabular} & \begin{tabular}[c]{@{}l@{}}300\\ 300\end{tabular} \\ \hline
\end{tabular}}
\end{table}

\subsection{Damped harmonic motion}

\begin{figure}[htbp]
\centering
\includegraphics[scale=0.3]{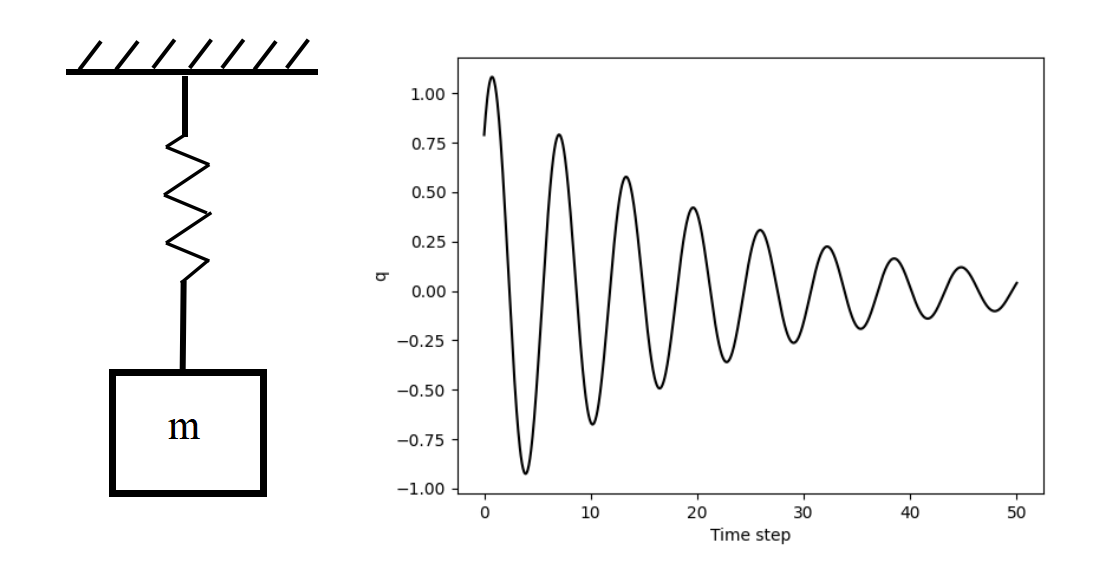}
\caption{Damped Harmonic Motion}
\label{damped}
\end{figure}

Perhaps the simplest example of a non-conservative system that can be described by the generalized Lagrange's equation is Damped Harmonic Motion, which is expressed as
$$\Ddot{q} + a\Dot{q} + k^{2}q = 0,$$
where $a$ represents the coefficient of friction, and $k$ denotes the coefficient of elasticity. In Figure \ref{damped}, we illustrate the model of this system and the variation of its coordinate $q$ over time.

Despite several classical expressions of generalized Lagrange's equations being proposed for Damped Harmonic Motion, we will provide an expression for the generalized Lagrange's equation based on Lemma \ref{generalL} introduced earlier, ensuring that $\mathscr{L}$ exactly represent the difference between the system's kinetic and potential energies.

In this example, the kinetic energy is $\dfrac{1}{2}m\Dot{q}^{2}$, and the potential energy is represented by elastic potential energy $\dfrac{1}{2}k^{2}q^{2}$. Hence, we have $$\mathscr{L} = \frac{1}{2}m\Dot{q}^{2} - \frac{1}{2}k^{2}q^{2}$$, leading to the equation $$F = \frac{d}{dt}\left( \frac{\partial \mathscr{L}}{\partial \Dot{q}}\right) - \frac{\partial \mathscr{L}}{\partial q} = m\Ddot{q} - kq.$$

In the numerical experiment, we choose the initial parameter values as a = 0.02 and k = 1. To generate training data, we select 40 trajectories, each starting from a randomly chosen point within the range of $[-1,1]^{2}$. For each trajectory, we sample 200 points with a step size of $h = 0.05$, resulting in a time interval of $T = 10$. These 8000 point pairs $(x_{t}, x_{t+1})$ are used as training data. Subsequently, we randomly split the training data into a training set and a testing set in a 1:1 ratio.

We use two networks to learn this system: GLNNs and the Baseline model (Baseline 2 model mentioned in Section 4.1). For the Baseline model, we utilize a three-layer fully connected network with a hidden dimension of 200. As for the GLNNs model, we employ two three-layer fully connected networks with a hidden dimension of 200. Both networks are trained using an Adam optimizer with a batch size of 1000 and a learning rate of 0.001. 

\begin{figure}[htbp]
\centering
\includegraphics[width=1.0\textwidth]{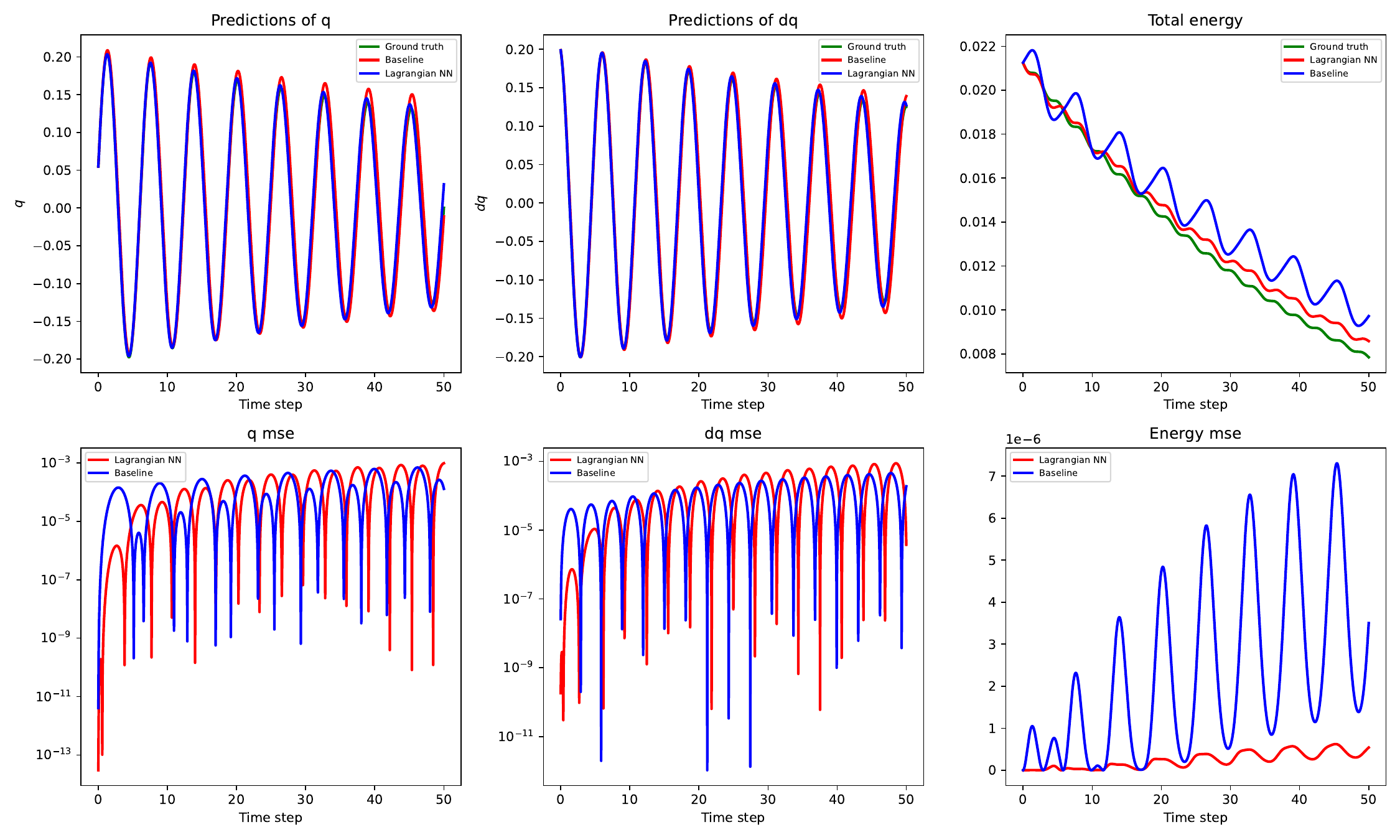}
\caption{Prediction of Damped Harmonic Motion}
\label{predamped}
\end{figure}

In Figure \ref{predamped}, we present the prediction results for Damped Harmonic Motion. The first row contains the first two plots showing the position and acceleration predictions for $q$. In the short-term prediction at $T = 50$, neither network demonstrates a significant superiority. However, the third plot in the first row displays the prediction for the total energy of the system. In contrast to the position in the phase flow, GLNNs starts to exhibit a noticeable advantage in energy prediction. In the second row, we plot the mean squared error (MSE) between the predicted energy of both networks and the ground truth. Through the MSE plot, we can clearly observe the higher accuracy of GLNNs in energy prediction.

\subsection{Compound double pendulum with friction}

\begin{figure}[htbp]
\centering
\includegraphics[width=0.8\textwidth]{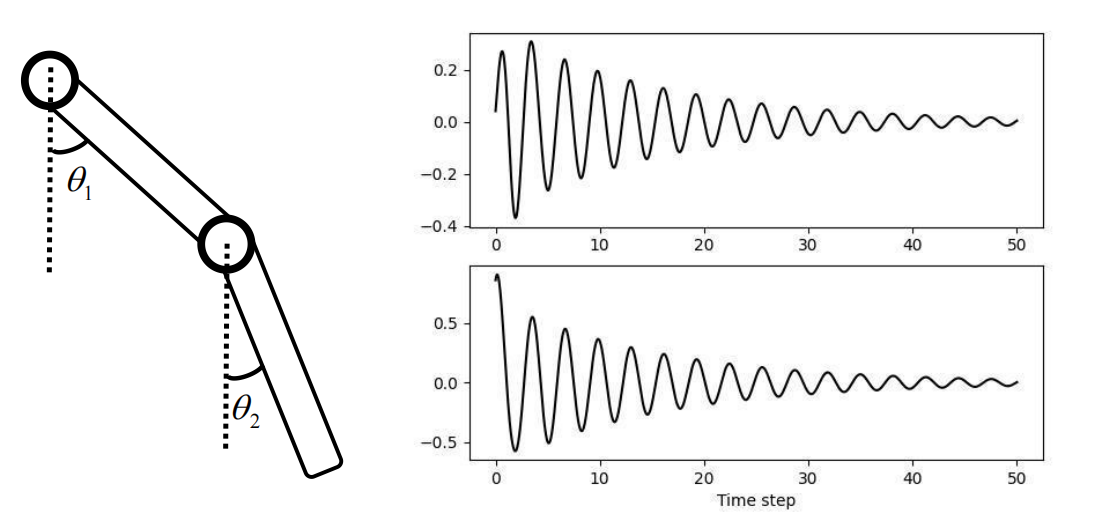}
\caption{Double pendulum with friction}
\label{dlp}
\end{figure}

Another non-conservative system example we consider is the double pendulum model. Unlike the conventional conservative double pendulum, we adopt the compound double pendulum with friction model proposed by Williams\cite{williams2023compound}. In this model, the two mass points are replaced with irregular-shaped rigid bodies, and frictional forces are considered within the system.

Next, we'll introduce the system. Consider a double pendulum with two identical rigid bodies of mass $m$ and irregular shape (see Fig \ref{dlp}). Note that the bodies need not be rods or squares. The upper body rotates around a fixed axle at $P$ in the x-direction, while the lower body rotates about an axle at $Q$ passing through the upper body. The massless axles are positioned at the same points within their respective plate's geometries. The bodies have uniform mass density distributions and equal moments of inertia about their axles. In our generalized setting, we account for friction and other nonconservative forces at the axles.

\begin{table}[htbp]
\caption{Parameters of the system}
\centering
\begin{tabular}{|l|l|l|}
\hline
Quantity               & Dimension                   & Description                                             \\ \hline
$m$                      & kg                          & Mass of first and second pendulum                       \\
$\theta_{1}$ & Dimensionless               & Angle which first pendulum makes with y axis            \\
$\theta_{2}$ & Dimensionless               & Angle which second pendulum makes with y axis           \\
$d$                      & m                           & Distance between the two axles                          \\
$I$                      & kg $m^{2}$     & Moment of inertia                                       \\
$c$                      & m                           & Distance from an axle to the centre of mass of pendulum \\
$\gamma_{1}$ & Dimensionless               & Friction coefficient for first axle                     \\
$\gamma_{2}$ & Dimensionless               & Friction coefficient for second axle                    \\
$g$                      & m $s^{-2}$ & Acceleration due to gravity                             \\ \hline
\end{tabular}
\end{table}

Here, we omit the specific derivation process of the system's equations of motion. Referring to the results presented in the article by Williams, we can obtain the equations of motion for the system as follows:
\begin{equation}
\begin{split}
     -\gamma_{1}\Dot{\theta_{1}} - \gamma_{2}(\Dot{\theta_{1}} - \Dot{\theta_{2}}) & = 2md^{2}\Ddot{\theta_{1}} + I\Ddot{\theta_{1}} + md^{2}\Ddot{\theta_{2}}cos(\theta_{1} - \theta_{2}) + md^{2}\Dot{\theta_{2}^{2}}sin(\theta_{1} - \theta_{2}) + 2mgdsin\theta_{1},\\
     \gamma_{2}(\Dot{\theta_{1}} - \Dot{\theta_{2}}) & = md^{2}\Ddot{\theta_{2}} + I\Ddot{\theta_{2}} + md^{2}\Ddot{\theta_{1}}cos(\theta_{1} - \theta_{2}) - md^{2}\Dot{\theta_{1}^{2}}sin(\theta_{1} - \theta_{2}) + mgdsin \theta_{2}.
\end{split}
\end{equation}
By solving the aforementioned equations involving $\Ddot{\theta_{1}}$ and $\Ddot{\theta_{2}}$, we can obtain an expression in the form of 
$$\Ddot{\theta_{1}} = F_{1}(\theta_{1},\theta_{2},\Dot{\theta_{1}},\Dot{\theta_{2}}), \Ddot{\theta_{2}} = F_{2}(\theta_{1},\theta_{2},\Dot{\theta_{1}},\Dot{\theta_{2}}).$$

Next, using the approach outlined in Lemma \ref{generalL}, we provide an expression for a generalized Lagrange's equation for this system. In this system, the potential energies of the double pendulum are denoted as 
\begin{equation}
\begin{split}
&T_{1} = \frac{1}{2}m(\Dot{x_{1}} + \Dot{y_{1}})^{2} + \frac{1}{2}J_{1}\Dot{\theta_{1}}^{2}\\
&T_{2} = \frac{1}{2}m(\Dot{x_{2}} + \Dot{y_{2}})^{2} + \frac{1}{2}J_{1}\Dot{\theta_{2}}^{2}
\end{split}
\end{equation}
where $J_{1} = \dfrac{1}{3}m\left( \dfrac{c}{2} \right)^{2}$ and $J_{2} = \dfrac{1}{3}m\left( \dfrac{c}{2} \right)^{2}$ represent the moments of inertia, and $(x_{1}, y_{1}), (x_{2}, y_{2})$ represent the traditional coordinates of the double pendulum. The potential energies of the double pendulum are represented as $U_{1} = mgy_{1}$ and $U_{2} = mgy_{2}$ . Thus, we have 
\begin{equation}
\begin{split}
     \mathscr{L} & = T-U  = \frac{1}{2}m(\Dot{x_{1}} +  \Dot{y_{1}})^{2} + \frac{1}{2}m(\Dot{x_{2}} + \Dot{y_{2}})^{2} + \frac{1}{6}m\left( \dfrac{c}{2} \right)^{2}(\Dot{\theta_{1}}^{2} + \Dot{\theta_{2}}^{2}) - mg(y_{1} + y_{2}),\\
    F_{1} & =  \frac{d}{dt}\left( \frac{\partial \mathscr{L}}{\partial \Dot{\theta_{1}}}\right) - \frac{\partial \mathscr{L}}{\partial \theta_{1}} \\
     &= \frac{d}{dt}\left(\frac{m^{2}c^{2}}{4}cos^{2}\theta_{1}\Dot{\theta_{1}} + \frac{m^{2}c^{2}}{4}cos\theta_{1}sin\theta_{1}\Dot{\theta_{1}}+\frac{mc^{2}}{12}\Dot{\theta_{1}} \right) - \frac{mc}{2}g sin\theta_{1} + \frac{mc^{2}}{4}(cos\theta_{1}sin\theta_{1}\Dot{\theta_{1}} - cos\theta_{1}\Dot{\theta_{1}}^{2})\\
     &= \frac{m^{2}c^{2}}{4}\Dot{\theta_{1}}^{2}(cos^{2}\theta_{1} - \sin^{2}\theta_{1} - 2sin\theta_{1}cos\theta_{1}) + \frac{mc^{2}}{4}(cos\theta_{1}sin\theta_{1}\Dot{\theta_{1}} - cos\theta_{1}\Dot{\theta_{1}}^{2}) - \frac{mc}{2}g sin\theta_{1},\\
    F_{2} & =  \frac{d}{dt}\left( \frac{\partial \mathscr{L}}{\partial \Dot{\theta_{2}}}\right) - \frac{\partial \mathscr{L}}{\partial \theta_{2}} \\
     &= \frac{d}{dt}\left(\frac{m^{2}c^{2}}{4}cos^{2}\theta_{2}\Dot{\theta_{2}} + \frac{m^{2}c^{2}}{4}cos\theta_{2}sin\theta_{2}\Dot{\theta_{2}}+\frac{mc^{2}}{12}\Dot{\theta_{2}} \right) - \frac{mc}{2}g sin\theta_{2} + \frac{mc^{2}}{4}(cos\theta_{2}sin\theta_{2}\Dot{\theta_{2}} - cos\theta_{2}\Dot{\theta_{2}}^{2})\\
     &= \frac{m^{2}c^{2}}{4}\Dot{\theta_{2}}^{2}(cos^{2}\theta_{2} - \sin^{2}\theta_{2} - 2sin\theta_{2}cos\theta_{2}) + \frac{mc^{2}}{4}(cos\theta_{2}sin\theta_{2}\Dot{\theta_{2}} - cos\theta_{2}\Dot{\theta_{2}}^{2}) - \frac{mc}{2}g sin\theta_{2}.
\end{split}
\end{equation}

Next, we proceed with numerical experiments, selecting initial parameters as $m = 1$, $c = 1$, $d = 1$, $g = 10$, $\gamma_{1} = 0.5$, and $\gamma_{2} = 0.5$. To generate training data, we choose 20 trajectories, each starting from a randomly selected point $(\theta_{1}^{0},\theta_{2}^{0},\Dot{\theta_{1}^{0}}, \Dot{\theta_{2}^{0}})$ within the range of $[1, -1]^{4}$. For each trajectory, we sample 500 points with a step size of $h = 0.02$, resulting in a time interval of $T = 10$. These 10,000 point pairs $(x_{t}, x_{t+1})$ are used as training data. Subsequently, we randomly split the training data into a training set and a testing set in a 1:1 ratio.

\begin{figure}[htbp]
\centering
\includegraphics[width=1.0\textwidth]{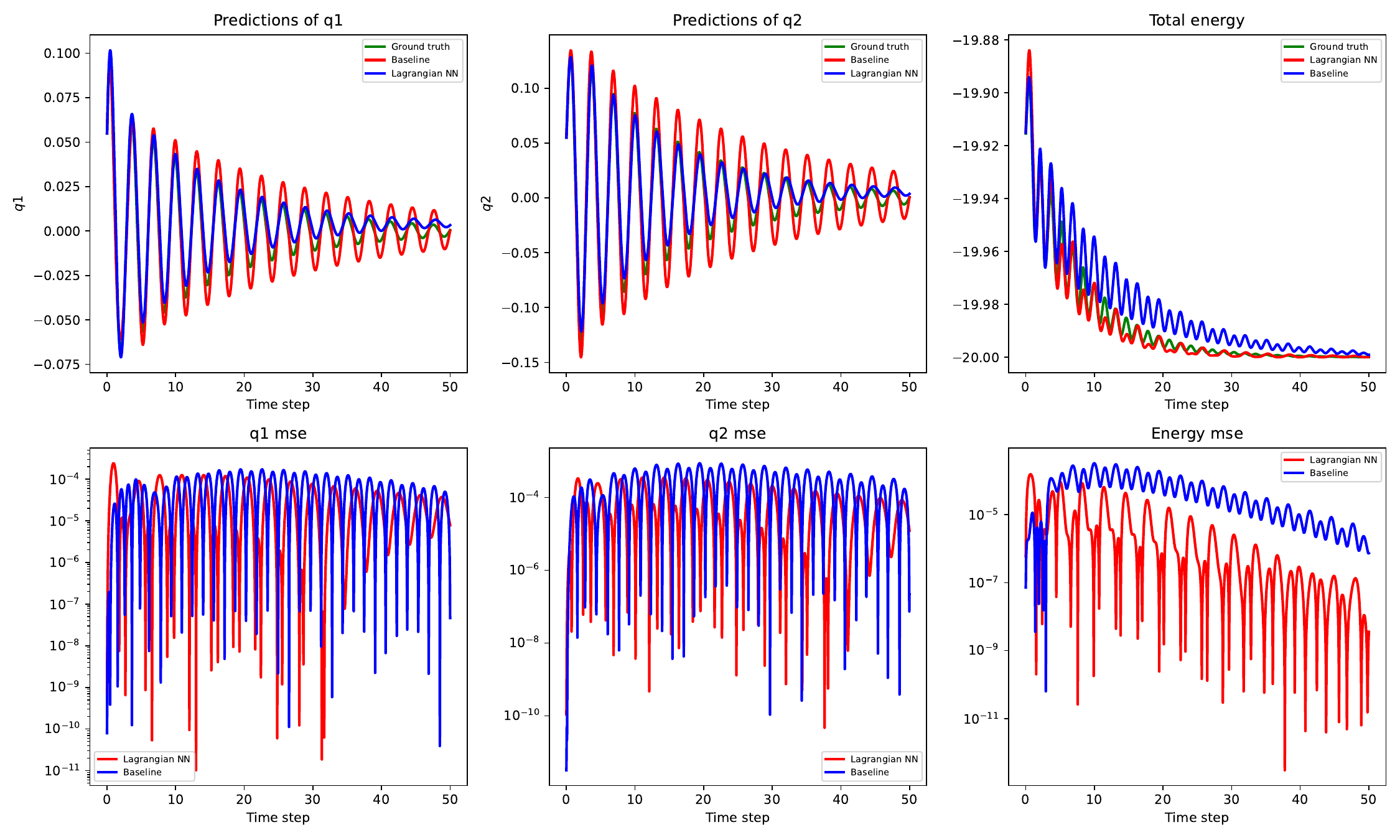}
\caption{Double pendulum with friction}
\label{predlp}
\end{figure}

Also, we test two network models. In the case of the Baseline model, a three-layer fully connected network with a hidden dimension of 200 is utilized. In the GLNNs model, two four-layer fully connected networks with a hidden dimension of 200 are employed. Both networks are trained using an Adam optimizer with a batch size of 1000 and a learning rate of 0.001.

In Figure 5, we present the prediction results for the compound double pendulum with friction. The first two plots in one row depict the angle predictions for $q_{1}$ and $q_{2}$ of the pendulum. The third plot shows the prediction for the system's energy. It is evident that in the case of the compound double pendulum, GLNNs also performs better in energy prediction.

\subsection{Investigation into the hyper-parameters}

In this subsection, we conduct experiments to investigate influence the hyper-parameters of our networks of predictions. 

We main consider change the number of layers and the size of hidden layers, we test layers of 2,3,4,5 and test the size of hidden layers of 50,100,200,400.

\begin{figure}[htbp]
\centering
 \includegraphics[width=1.0\textwidth]{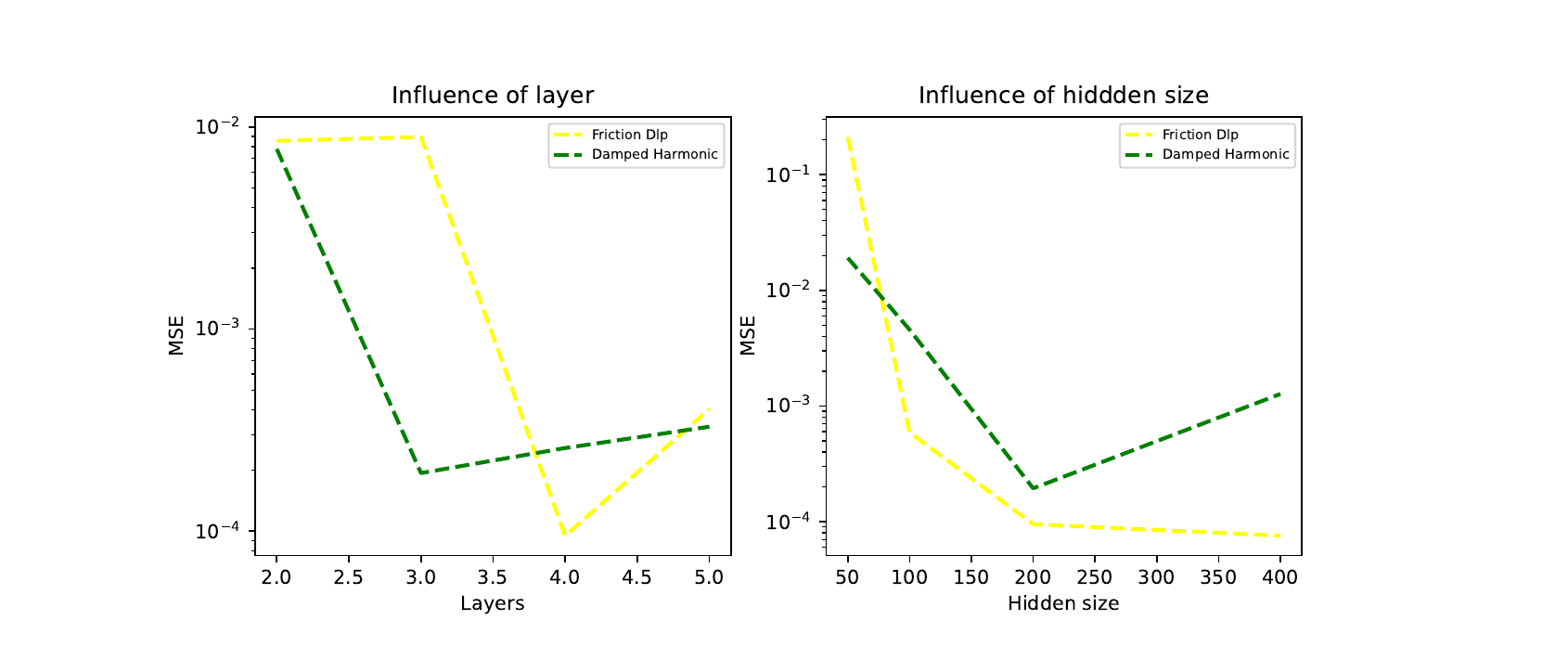}
\caption{Influence of hyper-parameters}
\label{para}
\end{figure}

In Figure \eqref{para}, we present the experimental results. It is evident that for the Damped Harmonic system, a network depth of two layers fails to achieve satisfactory training outcomes. However, increasing the depth to 4 or 5 layers with our available dataset does not lead to optimal training. Consequently, for the Damped Harmonic system, we opt for a three-layered GLNNs architecture, supplemented by a hidden size of 200. For the Compound Double Pendulum with friction system, given the increased dimensionality of the data, optimal training performance is observed at a network depth of four layers. While a hidden size of 400 offers marginal improvements over 200, the associated increase in computational cost leads us to maintain a hidden size of 200.

\begin{table}[htbp]
\centering\
\caption{Influence of hidden size}
\begin{tabular}{|l|l|l|l|l|}
\hline
Hidden size    & 50 & 100 & 200 & 400 \\ \hline
Damped harmonic (layer = 3) & 1.91e-2  & 4.57e-3   & 1.94e-4  & 1.27e-3   \\ \hline
Compound double pendulum (layer = 4) & 2.12e-1  & 5.91e-4   & 9.50e-5   & 7.55e-5   \\ \hline
\end{tabular}
\end{table}

\begin{table}[htbp]
\centering
\caption{Influence of layers}
\begin{tabular}{|l|l|l|l|l|}
\hline
Layer(hidden size = 200)   & 2 & 3 & 4 & 5 \\  \hline
Damped harmonic          & 7.79e-2  & 1.94e-4   & 2.58e-4 & 3.29e-4   \\ \hline
Compound double pendulum & 8.55e-3  & 8.93e-3  & 9.50e-5 & 4.02e-4  \\ \hline
\end{tabular}
\end{table}

\section{Conclusions}

We outline the primary contributions of our research as follows:

\begin{itemize}
\item We introduce a methodology to transform physically motivated dissipative systems into a generalized Lagrangian framework.
\item We develop Generalized Neural Networks (GLNNs) tailored for physically motivated dissipative systems.
\item We derive the Lagrangian representation for a compound double pendulum subject to friction.
\item Through various experiments, we demonstrate the superior performance of GLNNs in both one-dimensional and two-dimensional scenarios.
\end{itemize}

Based on our theoretical findings and experimental results, we ascertain the following advantages associated with GLNNs:

\begin{itemize}
\item In contrast to LNNs, GLNNs are versatile, applicable to both conservative and non-conservative systems.
\item When addressing real-world physical systems with frictional components, Generalized Lagrangian Neural Networks (GLNNs) offer improved accuracy in predicting system energy and phase flow dynamics.
\item The inherent preservation of the system's Lagrangian structure in GLNNs enhances their predictive efficacy in specific contexts.
\end{itemize}

While GLNNs offer certain advantages, it is crucial to acknowledge their inherent limitations. One notable limitation is that, when applied to systems characterized by energy dissipation, GLNNs cannot accurately ensure that the predicted energy diminishes entirely. The forecasting process may exhibit instances of energy rebound, deviating from expected physical characteristics. Another constraint lies in the elevated training complexity associated with GLNNs compared to baseline models. Given the increased network complexity, there is an inevitable escalation in both training time requirements and computational resources.

\hspace*{\fill}\\
\textbf{\textbf{Acknowledgments}}

This research is supported by National Natural Science Foundation of China (Grant Nos. 12171466 and 12271025).

\end{CJK}

\end{document}